\documentclass[a4paper,12pt]{article}
\usepackage{amssymb}
\usepackage{amsmath}
\usepackage{color}
 \usepackage{amsthm}
 \usepackage{CJK}
 \usepackage{geometry}
\usepackage{graphics}
 \usepackage[dvips]{graphicx}
\usepackage{multirow}
\usepackage{fancyhdr}
\usepackage{multirow}
\usepackage{booktabs}
\usepackage{longtable}
\usepackage{diagbox}
\usepackage{ltablex}
\usepackage{subcaption}
\usepackage{geometry}

\newtheorem{lem}{Lemma}[section]
\newtheorem{thm}[lem]{Theorem}

\newtheorem{cor}[lem]{Corollary}

\newtheorem{defi}[lem]{Definition}

\newtheorem{problem}[lem]{Problem}
\newtheorem{proposition}[lem]{Proposition}

\newcolumntype{C}{>{\centering\arraybackslash}X}

\begin{document}
	
\title{Probability model of connectedness for graphs with logarithmical degree}
\author{Huanshen Jia$^1$ and Jianguo Qian$^{1,2}$\footnote{Corresponding author. E-mail:jgqian@xmu.edu.cn (J.G. Qian)}\\
	 \small $^1$School of Mathematics and Statistics, Qinghai Minzu University, P. R. China\\
    \small $^2$School of Mathematical Sciences, Xiamen University, P. R. China}
	\date{}
	\maketitle
\begin{abstract}
For a graph $G$ of order $N$, an integer $k$ with $0\leq k< N$ and a real number $p$ with $0\leq p\leq 1$, let  $p(G,k)$ and  $t(G,p)$ be the probabilities that $G$ is disconnected when $k$ vertices are removed randomly and when each vertex is removed independently with probability $p$, respectively.  We deduce a formula of $p(G,k)$ for Harary graphs and show that the phase transition on $p(G,k)$ and $t(G,p)$ occurs in the class that has logarithmical  degree. We also perform simulations on $p(G,k)$ for the Hypercube and randomly generated graphs of the same scale, which suggests that the regularity is a key structural feature that guarantees nice performance on $p(G,k)$.  In particular, the  Hypercubes have better performance than regular graphs, on average. Surprisingly, Harary graphs, known to be super-connected, have neither better performance than  Hypercubes nor even than the average of regular graphs.
 \end{abstract}
\noindent{\bf Keywords:} Vertex-failure;  Graph connectedness; Logarithmical degree;  Probability
\section{Introduction}
In graph theory, many types of measurements concerning graph connectedness were introduced from both theoretical and practical sides.  One of the frequently used is the traditional {\it vertex-connectivity} or simply {\it connectivity}, which is defined as the maximum number of vertices whose removal does not  disconnect  the graph. By Menger's theory, if a graph has vertex-connectivity $k$, then any two vertices are connected by at least $k$  internally vertex-disjoint paths \cite{Menger} and, vice versa.

The vertex-connectivity ensures only the maximum number of the vertices whose removal guarantees the connection of a graph. However, even if more vertices are removed, a graph may still remain connected. It is then interesting to quantify the effect of vertex removal or, in the sense of probability model, to consider the problem: what is the probability that a graph is connected (or disconnected) when a given number of vertices are removed? Further, if each vertex has a certain failure probability independent of one another, what is the probability that the graph is connected? However, the number of vertex subsets in a graph increases exponentially when the scale of the graph increases. Therefore, the two problems are even more difficult for a large-scale graph in general.

In practice, in order to avoid the combinatorial difficulty of the two problems, one might resort to some numerical simulation such as Monte Carlo sampling techniques, finding efficient algorithms for special classes of graphs, or using bounding techniques \cite{Weichenberg}. Another direction is to consider the problems from the perspective of random graph \cite{Erdos,Nikoletseas,Takabe}. In application, this model has also attracted particular interest  in large scale wireless communication networks \cite{Kawabata,Kuo,Nozaki,Sudevalayam,Yick}.

For a graph $G$ with $N$ vertices and positive integer $k$ with $0\leq k< N$, we denote by  $p(G,k)$ the probability that $G$ is disconnected when $k$ vertices are removed randomly.  The measurement  $p(G,k)$ has close relation with the number of the connected sets in a graph, which was extended from the well-known subtree-counting problem and has  received much attention in recent ten years \cite{Cambie,Kangas,Vince}. A vertex subset of a graph $G$ is called a  {\it connected set} if it induces a connected graph. Therefore,  the number of the vertex subsets of order $k$ in a graph of order $N$ whose removal remains a connected graph equals the number of the connected sets of order $N-k$. However, except the few graphs e.g., paths, cycles, complete (bipartite) graphs,  wheels, necklaces, ladders, prism \cite{Vince}  and criss-cross prism \cite{Haslegrave}, this problem still remains open even for very special graphs such as the Hypercubes and the grid graphs (the Cartesian product of two paths) \cite{Vince}. For general graphs, in \cite{Vince} Vince posed a question whether the number of the connected sets is necessarily unimodal? Several lower and upper bounds on the the maximum number of connected sets in $d$-regular graphs were also established by asymptotical analyses, computer search and constructing particular graphs \cite{Cambie,Kangas,Vince}.

 For a real number $p$ with $0\leq p\leq 1$, we denote by $t(G,p)$ the probability that $G$ is disconnected when each vertex is removed independently with probability $p$ (also called the {\it node breakdown probability \cite{Takabe}} in networking analysis). We note that the probabilities  $t(G,p)$ and  nRel$(G; p)=1-t(G,1-p)$ are also called the {\it network breakdown probability} \cite{Takabe,Nozaki} and {\it node reliability} \cite{Matthew}, respectively. It is well known that the function  nRel$(G; p)$ in $p$ is not necessarily monotone increasing on the interval $(0,1)$ in general, which is very different from its edge version, namely the {\it all-terminal reliability} defined as the probability that $G$ is connected when each edge is removed independently with probability $1-p$ which was initially introduced for optical networks with passive nodes and vulnerable fiber under stress of disconnection  \cite{Weichenberg}. In terms of  $t(G,p)$, it was shown \cite{Brown} that $t(G,p)$ has a decreasing interval in $(0,1)$ for sufficiently sparse graphs, e.g., trees, and is concave up when $p$ is sufficiently close to 1 for any graph and, moreover, $t(G,p)$ has at least two fixed points (the point $p$ satisfying  $t(G,p)=p$) in $(0,1)$ for any sufficiently large graph with bounded maximum degree.  In \cite{Nozaki}, an upper bound on $t(G,p)$ was given, which is very close to the experimental values of $t(G,p)$ when $p$ is small.

 In the topic concerning the graph connectedness, semi-dense graphs, that is, graphs with logarithmical average degree with respect to the number of vertices, plays an important role, in which many  phase transition phenomena were found. The first work may date back to 1959 \cite{Erdos} on Erd\H{o}s-R\'{e}nyi  random graphs (ER model), which showed that the logarithmical average degree is the threshold for random graphs to be connected, namely the so-called percolation transition \cite{Takabe}. Specifically, for a random graph $G$ with $N$ vertices and $M$ edges, $G$ is connected if $M> (N/2)\ln N$ and disconnected if $M<(N/2)\ln N$ asymptotically almost surely, as $N\rightarrow+\infty$. Based on the correspondence between configurations and random graphs provided by
Bollob\'{a}s, in \cite{Nikoletseas}  Nikoletseas et al. revealed the threshold on the connectivity of the random regular semi-dense graphs with edge removing.  In \cite{Takabe}, using mean-field approximation, Takabe et al. obtained an approximation formula of $t(G,p)$ for regular random graphs, which predicts that the phase transition for $t(G,p)$ occurs at $p=e^{-1/d}$ in $d\ln N$-regular random graphs as $N\rightarrow+\infty$. The graphs with logarithmical average degree are also of interest in the complex network analysis, because probabilistic analyses show
the existence of phase transition phenomena of some crucial problems such as community detection \cite{Abbe,Takabe}. Even so, we are not aware of the explicit formula of  $p(G,k)$ and $t(G,p)$ for any graph family with logarithmical average degree, up to now.

 Due to their relatively small link complexity, the graphs with logarithmical average degree also extensively applied in the design of communication networks. The most typical example would be the Hypercube and its variants. In fact, in addition to its lower degree,  Hypercubes have also many other attractive properties such as simple recursive structure, simple labeling scheme, small diameter, high symmetry, routing capabilities \cite{Beivide,Boesch1,Day,Jia,Li2,Qiao,Vaidya} and nice structure mapping property \cite{Saad}. Another example we focus on is the Harary graph, which has also been shown to possess attractive reliability properties and therefore considered as candidate
topologies for high-reliability applications. It is well known that the vertex and edge connectivity of a graph are both bounded from above by the average degree of the graph and can be achieved through the construction of Harary graphs \cite{Weichenberg}. Moreover, for edge connectivity $\lambda$, it was shown \cite{Deng} that the Harary graph has minimum values on the edge version of $p(G,k)$ for $k=\lambda,\lambda+1,\ldots,2\lambda-3$ over all super-$\lambda$ connected graphs (every minimum edge cut isolates a vertex of degree $\lambda$) of the same scale. 

 In this paper, we consider the regular graphs with  logarithmical degree. More specifically, we focus on the class of $n$-regular graphs of order $2^n$. For comparison, we also choose Harary graph and Hypercube of the same scale as two concrete instances in the class. In the following section, we deduce an explicit formula of $p(G,k)$ for Harary graphs. Further, by asymptotical analyses we show that the phase transitions on $p(G,k)$ and $t(G,p)$ occurs at  $k=2^{n-2}$ and $p=0.25$ as $n\rightarrow \infty$ for $n$-regular  Harary graph of order $2^n$, respectively. To our knowledge up to now, this would be the first known infinite graph class with  logarithmical degree that has such phenomenon. In the third section, we perform simulations of $p(G,k)$ for the Hypercube of order $2^6$ and 50 randomly generated graphs and regular graphs of the same scale. Our numerical results indicate that the regularity is apparently a key structural feature that guarantees lower value of  $p(G,k)$.  In particular, in addition to their well-structured properties for network application,  the  Hypercube has also lower values of $p(G,k)$ than the average of the 50 randomly generated regular graphs.  Although Harary graphs are super connected and have minimum values on the edge version of $p(G,k)$ for some $k$, it is a little surprising that Harary graphs have neither lower value on $p(G,k)$ than Hypercubes does nor, even than the regular graphs in average.

 \section{Theoretical results on Harary graphs}

In this section we consider Harary graphs, which might provide us an intuitive insight into the two measurements $p(G,k)$ and $t(G,p)$.
\begin{defi}\label{harary}{\rm({\bf Harary Graph} \cite{Boesch2,Weichenberg})}
For positive integers $N$ and $w$ with $1\leq w\leq N/2$, the Harary Graph $H(N,w)$ is a $2w$-regular graph consisting of $N$ vertices labelled by $1,2,\cdots,N$ such that the vertex $i$ is adjacent to $j$ if and only if $j\in\{i\pm 1,i\pm 2,\ldots,i\pm w\}\ ({\rm mod}\ N)$.
\end{defi}

We note that  Harary graphs are particular circulant graphs and $2w$-regular with vertex-connectivity $2w$. 

For a vertex $v\in V(G)$,  the {\it degree} of $v$ is defined as usual as the number of the vertices adjacent to $v$. In particular, if every vertex has the same degree $k$ then we call $G$ $k$-{\it regular}. A subset $W$ of $V(G)$ is called a {\it vertex cut} of $G$ if $G-W$ is not connected. 
For a polynomial $f(x)$ and an integer $k$, we denote by $[x^k]f(x)$ the coefficient of $x^k$ in $f(x)$.
\begin{proposition}\label{coef} For  any positive integers $k,s$ and $t$, we have
 $$[x^k]\left(\sum\limits_{i=0}^sx^i\right)^t=\sum\limits_{j=0}^{\lfloor k/(s+1)\rfloor}(-1)^j{t\choose j}{k-j(s+1)+t-1\choose t-1}.$$
\end{proposition}
\begin{proof}
 $$\left(\sum\limits_{i=0}^sx^i\right)^t=(1-x^{s+1})^t(1-x)^{-t}=\sum\limits_{i=0}^t(-1)^i{t\choose i}x^{i(s+1)}\sum\limits_{j=0}^{+\infty }{j+t-1\choose t-1}x^j,$$
where the last equality holds by the generalized binomial theorem. Hence,
$$[x^k]\left(\sum\limits_{i=0}^sx^i\right)^t=\sum\limits_{j=0}^{\lfloor k/(s+1)\rfloor}(-1)^j{t\choose j}{k-j(s+1)+t-1\choose t-1},$$
where $\lfloor k/(s+1)\rfloor$ in the right summation is because $k-j(s+1)+t-1\geq t-1$.
\end{proof}
\begin{thm}\label{hg} Let $N$ and $w$ be two positive integers with $w\leq N/2$. If $k<2w$ then $p(H(N,w),k)=0$, and if $2w\leq k<N$ then
\begin{equation}\label{hqe}
p(H(N,w),k)=1-\left(\alpha+\sum\limits_{i=1}^{w-1}i\alpha_i+N\sum\limits_{j=w}^k\alpha_j\right)/{N\choose k},
  \end{equation}
where
$$\alpha=\sum\limits_{j=0}^{\lfloor k/w\rfloor}(-1)^j{N-k\choose j}{N-jw-1\choose N-k-1}$$
and
$$\alpha_i=\sum\limits_{j=0}^{\lfloor (k-i)/w\rfloor}(-1)^j{N-k-1\choose j}{N-jw-i-2\choose N-k-2}.$$
 \end{thm}
\begin{proof} Since any vertex subset of $H(N,w)$ with fewer than $2w$ vertices is not a vertex cut, $p(H(N,w),k)=0$ follows directly if $k<2w$.

We now assume $k\geq 2w$. For any integer $i$ with $i\geq w$, we call $i$ vertices on $H(N,w)$ a {\it section} if they  are consecutive on  $H(N,w)$. By the definition of Harary graph, one can see that a vertex subset $K$ of  $H(N,w)$ is a vertex cut if and only if $K$ contains at least two disjoint sections $V_i$ and $V_j$, such that $N\setminus K$ has at least one vertex in each of the two parts of $H(N,w)$ separated by  $V_i$ and $V_j$.

Let $K$ be a vertex subset of order $k$ in $H(N,w)$ but not a vertex cut, that is, $V(H(N,w))\setminus K$ is a connected set. If $K$ does not contain any section and $1\notin K$, then the choice of $K$ is equivalent to inserting at most $w-1$ vertices into each of the $N-k$ intervals on the cycle separated by the vertices of $V(H(N,w))\setminus K$ (label an arbitrary vertex as the vertex 1). The number of such choices, i.e., $\alpha$, equals the coefficient of $x^k$ in the following generating function:
\begin{equation}\label{a}
(1+x+x^2+\cdots+x^{w-1})^{N-k}
\end{equation}

If $K$ does not contain any section but $1\in K$, then the vertex 1 is in an interval on the cycle separated by the vertices of $V(H(N,w))\setminus K$. Label this interval by $I$ and assume that $I$ contains $i$ vertices of $K$, where $1\leq i\leq w-1$. Such choice is equivalent to inserting at most $w-1$ vertices into each of the $N-k-1$ intervals (other than $I$) on the cycle separated by the vertices of  $V(H(N,w))\setminus K$, such that the total number of the vertices in the $N-k-1$ intervals equals $k-i$. The number of such choices, i.e., $\alpha_i$, equals the coefficient of $x^{k-i}$ in the following generating function:
\begin{equation}\label{ai}
(1+x+x^2+\cdots+x^{w-1})^{N-k-1}
\end{equation}
In addition, each vertex in the interval $I$ could be serve as vertex 1. Hence, the total number of such choices equals $i$ times the above coefficient.

We now assume that $K$ contains exactly one section $V_j$ and assume that  $V_j$ consists of $j$ vertices where $j\geq w$. In this case, except the interval $V_j$, the $N-k$ vertices of $V(H(N,w))\setminus K$ forms $N-k-1$ intervals. Insert the $k-j$ vertices of $K\setminus V_j$ to these $N-k-1$ intervals so that each interval is inserted at most $w-1$ vertices, the number of which is the coefficient of $x^{k-j}$ in (\ref{ai}). Further, $V_j$ has $N$ choices (each of $1,2,\ldots,N$ could be as the starting vertex of the section $V_j$).

In conclusion, the theorem follows by Proposition \ref{coef}.
 \end{proof}

In the following we consider the asymptotic behaviour  of $p(H(N,w),k)$ as $w\rightarrow\infty$ for the special case that $N=4^w$. The following theorem shows that  $p(H(N,w),k)$  has phase transition at $k=4^{w-1}$, which would be the first known example of the graph classes of logarithmical average degree that has such phenomenon.
\begin{thm}\label{asym}Let $k=c4^w$ where $0<c<1$. Then, as $w\rightarrow\infty$,
$$p(H(4^w,w),k)\sim 1-(1+(1-c)(4c)^w)e^{-(1-c)(4c)^w}$$
and consequently,
\begin{align*}p(H(4^w,w),k)\rightarrow
\begin{cases}
0, & {\rm if}\ 0<c<\frac{1}{4}, \\
1-\frac{7}{4}e^{-3/4}\approx 0.1733, & {\rm if}\  c=\frac{1}{4},\\
1, & {\rm if}\  \frac{1}{4}<c<1.
\end{cases}
\end{align*}
\end{thm}
\begin{proof}By Theorem \ref{hg}, we assume $k\geq 2w$. Consider firstly the  asymptotic value of $\alpha$, that is, the $k$-th derivative of
$$\left(\sum\limits_{i=0}^{w-1}x^i\right)^{4^w-k}$$
at $x=0$. By Cauchy's integral formula, we have
$$\alpha=\frac{1}{2\pi i}\oint_C \left(\frac{1-x^{w}}{1-x}\right)^{4^w-k}\frac{dx}{x^{k+1}},$$
where $C$ is a small closed curve  in complex plane containing the origin. Let $F(x)=\left(\frac{1-x^{w}}{1-x}\right)^{4^w-k}\frac{1}{x^{k}}$ and choose $C$ as the  cycle with radius $r$ where $0<r<1$. Then $x=re^{i\theta}$ and
$$\alpha=\frac{1}{2\pi}\int_{-\pi}^{\pi}F(re^{i\theta})d\theta.$$

  We apply saddle point method to calculate the above integral. Let
$$f(x)=\ln F(x)=(4^w-k)\ln(1-x^w)-(4^w-k)\ln(1-x)-k\ln x.$$
Since $k=c4^w$, we have
\begin{equation}\label{fd}
f'(x)=-\frac{(4^w-k)wx^{w-1}}{1-x^w}+\frac{4^w-k}{1-x}-\frac{k}{x}
=\frac{(1-c)4^w}{x}\left(-\frac{wx^{w}}{1-x^w}+\frac{x}{1-x}-\frac{c}{1-c}\right),
\end{equation}
where $f'(x)$ denotes the derivative of $f(x)$.
Setting $f'(x)=0$, we have
\begin{equation}\label{fd1}
-\frac{wx^{w}}{1-x^w}+\frac{x}{1-x}-\frac{c}{1-c}=0.
\end{equation}

Since $w\rightarrow\infty$ and $|x|<1$, we have $wx^{w}/(1-x^w)\rightarrow 0$ and hence, the saddle point is $x_0=c$. Choose the radius $r=c$. Then $x=ce^{i\theta}$, and the saddle point  $x_0=c$ corresponds to $\theta_0=0$.  Consider the second derivative of $f(\theta)$ at $\theta=0$. By (\ref{fd}) and (\ref{fd1}),
$$f''(\theta)=(1-c)4^w\left(\frac{w^2c^we^{iw\theta}}{(1-c^we^{iw\theta})^2}-\frac{ce^{i\theta}}{(1-ce^{i\theta})^2}\right).$$

When $\theta=0$, we have
$$\left(\frac{w^2c^we^{iw\theta}}{(1-c^we^{iw\theta})^2}\right)/\left(\frac{ce^{i\theta}}{(1-ce^{i\theta})^2}\right)\sim w^2c^w/c\rightarrow 0.$$
Hence,
$$f''(0)\sim-4^w(1-c)\frac{c}{(1-c)^2}=-\frac{c4^w}{1-c}.$$
Since $f'(0)=0$, so by Taylor expansion, for small $\theta$, we have
$$f(\theta)=f(0)-\frac{c4^w}{2(1-c)}\theta^2+O(4^w\theta^3).$$
 Therefore,
$$\alpha=\frac{1}{2\pi}\int_{-\pi}^{\pi}F(re^{i\theta})d\theta=\frac{1}{2\pi}\int_{-\pi}^{\pi}e^{f(\theta)}d\theta
\sim\frac{e^{f(0)}}{2\pi}\int_{-\infty}^{\infty}e^{-\frac{c4^w}{2(1-c)}\theta^2}d\theta=\frac{e^{f(0)}}{\sqrt{2\pi c4^w/(1-c)}},$$
where the last equality holds by Gaussian integral. Notice that
$$f(0)=\ln F(c)=4^w(1-c)\ln(1-c^w)-4^w(1-c)\ln(1-c)-c4^w\ln c.$$
By Stirling's formula $p!\sim\sqrt{2\pi p}\left(\frac{p}{e}\right)^p (p\rightarrow \infty)$, for $c\in\{1/4^w,2/4^w,\ldots,1\}$ we have
\begin{equation}\label{stirl}
{4^w\choose c4^w}\sim\frac{1}{\sqrt{2\pi c4^w(1-c)}}e^{-4^w((1-c)\ln(1-c)+c\ln c)}.
\end{equation}
Further, notice that $\ln(1-c^w)\sim -c^{w}$ as $w\rightarrow\infty$. Therefore,
$$\alpha\sim(1-c){4^w\choose c4^w}e^{-(1-c)(4c)^w}.$$

We now consider $\alpha_i$. By the same argument as for $\alpha$ above, we have
$$\alpha_i\sim(1-c)c^i\alpha\sim(1-c)^2c^i{4^w\choose c4^w}e^{-(1-c)(4c)^w}.$$

Write $\beta={4^w\choose c4^w}$. Then by Theorem \ref{hg}, we have
$$p(H(4^w,w),k)=1-\frac{1}{\beta}\left(\alpha+\sum\limits_{i=1}^{w-1}(i \alpha_i)+4^w\sum\limits_{i=w}^{c4^w}\alpha_i\right).$$

By the asymptotic values of $\alpha$ and $\alpha_i$ we obtained above, it is not difficult to see that
\begin{equation}\label{alpha}
\alpha+\sum\limits_{i=1}^{w-1}i\alpha_i\sim \beta e^{-(1-c)(4c)^w}\ \ \mbox{\rm and}\ \ 4^w\sum\limits_{i=w}^{c4^w}\alpha_i\sim \beta(1-c)(4c)^we^{-(1-c)(4c)^w}.
\end{equation}
Therefore,
\begin{equation}\label{pHk}
p(H(4^w,w),k)\sim 1-(1+(1-c)(4c)^w)e^{-(1-c)(4c)^w},
\end{equation}
meaning that $p(H(4^w,w),k)=0, 1-7e^{-3/4}/4$ and 1 if $c<1/4, c=1/4$ and $c>1/4$, respectively. This is consistent with the fact in Theorem \ref{hg} that $p(H(4^w,w),k)=0$ when $k<2w$, since $2w/4^w\rightarrow 0<1/4$. This completes the proof.
\end{proof}

 Applying Theorem \ref{hg}, in Figure 1 we illustrate the values of $p(H(4^w,w),k)$ for $w=3,4,5,6,7$ and $k=s\times4^{w-3}, s=1,2,\ldots,63$, which intuitively appears the  behavior of the phase transition at $k=4^{w}/4=4^{w-1}$, i.e, $s=16$ as proved in Theorem \ref{asym}.

\begin{figure}[!htbp]
\begin{center}
\includegraphics[height=8cm]{1}\\
{\bf Figure 1.}\,\ The  values of $p(H(4^w,w),k)$ for $w=3,4,5,6,7$ and $k=s\times4^{w-3}, s=1,2,\ldots,63$.
\end{center}
\end{figure}

In the following we give the asymptotical behaviour of $t(G,p)$ on $H(4^w,w)$.
\begin{cor}\label{tp}As $w\rightarrow +\infty$,
\begin{align*}t(H(4^w,w),p)\rightarrow
\begin{cases}
0, & {\rm if}\ 0<p<\frac{1}{4}, \\
1-\frac{7}{4}e^{-3/4}\approx 0.1733, & {\rm if}\  p=\frac{1}{4},\\
1, & {\rm if}\  \frac{1}{4}<p<1.
\end{cases}
\end{align*}
\end{cor}
\begin{proof} By the definition of $t(G,p)$, it is clear that
$$t(H(4^w,w),p)=\sum\limits_{k=0}^{4^w}{4^w\choose k}p(H(4^w,w),k)p^k(1-p)^{(4^w-k)}.$$
Let
\begin{equation}\label{phi}
\phi_w(x)=(1+(1-x)(4x)^w)e^{-(1-x)(4x)^w}.
\end{equation}
Then by (\ref{pHk}), we have
$$t(H(4^w,w),p)\sim1-\sum\limits_{k=0}^{4^w}{4^w\choose k}\phi_w\left(\frac{k}{4^w}\right)p^k(1-p)^{(4^w-k)}.$$
Let
 $Y\sim {\rm Binomial}(4^w,p)$ (Binomial distribution). Then by the Law of Large Numbers, $Y/4^w\rightarrow p$ almost certainly. Hence,
$$\sum\limits_{k=0}^{4^w}{4^w\choose k}\phi_w\left(\frac{k}{4^w}\right)p^k(1-p)^{(4^w-k)}=\sum\limits_{k=0}^{4^w}\phi_w\left(\frac{k}{4^w}\right)\mathbb{P}\left(Y=k\right)
=\mathbb{E}[\phi_w(Y/4^w)].$$
For $z\geq 0$, we have $(1+z)e^{-z}\leq 1+e^{-1}<2$ and, hence $0\leq \phi_w(Y/4^w)\leq 2$, i.e., $\phi_w(Y/4^w)$ is upper bounded. Therefor, $\phi_w(Y/4^w)\rightarrow \phi(p)$ almost certainly as $w\rightarrow\infty$. So by the Dominated Convergence Theorem, we have
$$\lim\limits_{w\rightarrow\infty}\mathbb{E}[\phi_w(Y/4^w)]=\mathbb{E}[\lim\limits_{w\rightarrow\infty}\phi_w(Y/4^w)]=\phi(p),$$
where $\phi(p)=\lim\limits_{w\rightarrow\infty}\phi_w(p)$. The theorem follows directly by (\ref{phi}).
\end{proof}

Again applying Theorem \ref{hg}, in Figure 2 we illustrate the values of $t(H(4^w,w),p)$ for $w=3,4,5,6$, and $p=0,0.025,0.05,0.075,\ldots, 0.975,1$, which intuitively appears the phase transition on $t(G,p)$ at $p=0.25$ as proved in Theorem \ref{asym}.
\begin{figure}[htbp]
\begin{center}
\includegraphics[height=8.5cm]{2}\\
{\bf Figure 2.}\,\ The  values of $t(H(4^w,w),p)$ for $w=3,4,5,6$ and $p=0,0.025,0.050,\ldots, 0.975,1$.
\end{center}
\end{figure}

In \cite{Takabe}, Takabe et al. gave an  approximation analysis on $t(G,p)$ for regular random graphs, which predicts that the phase transition on $t(G,p)$ occurs at $p=0.5$  in $n$-regular random graphs of order $2^n$ as $n\rightarrow+\infty$. Corollary \ref{tp} shows that the phase transition in the Harary graphs $H(4^w,w)$ occurs at $p=0.25$, which is smaller than 0.5. This implies that, if the prediction of Takabe et al. is true, then the Harary graphs do not perform better behaviour on $t(G,p)$ than regular graphs of the sam scale, in average.

\section{Simulation results on Hypercube and\\ random graphs}

For a positive integer $n$, the $n$-dimensional hypercube $Q_n$ is an $n$-regular graph with $2^n$ vertices labeled by $n$-bit binary strings $x_1x_2\cdots x_n\in\{0,1\}^n$, such that there is an edge between any two vertices if and only if their binary strings differ by one and only one bit. In this section, we perform our simulation on $p(G,k)$ for Hypercube $Q_6$ and, in comparison, for 50 randomly generated  6-regular graphs with 64 vertices and 50 randomly generated graphs with the same number of vertices and edges (by Matlab PRNG, Pseudo-Random Number Generator based on Mersenne Twister Algorithm).

In our simulation, for each graph $G$ and each integer $k$  with $3\leq k\leq63$,  we choose 4000 graphs obtained from $G$ by randomly removing $k$ vertices. When $k=1$ and $k=2$, we choose every vertex and every pair of vertices as vertex faulty sets, respectively. Further, we denote by $\lambda_k(G)$, or simply $\lambda_k$, the number of disconnected graphs among all the 4000 faulty graphs with $k$ faulty vertices. It is clear that $p(G,k)\approx\lambda_k(G)/4000$. We note that, two vertices $u$ and $v$ in a graph $G$ of order $n$ are connected if and only if the element $(u,v)$ (the element at the $u$-th row and $v$-th column) in the matrix $A^{n-1}+A^{n-2}$ is not 0, where $A$ is the adjacency matrix of $G$.

We list the numerical results of $\lambda_k/4000$ for $Q_6$ and the average of 50 randomly generated graphs, denoted by $\overline{R}_6$ , and the  average of 50 randomly generated regular graphs, denoted by $\overline{R}R_6$, in Appendix (Table 1).

For an intuitive illustration, we plot the numerical values  of $\lambda_k/4000$ for $Q_6,\overline{R}_6$ and $\overline{R}R_6$, and the theoretical values of $p(G,k)$ for $H(4^w,w)$ we obtained in the previous section in Figure 3.

From Figure 3, we can see that the regularity is apparently a key structural feature that guarantees low value of  $p(G,k)$ in average. In particular, the  Hypercubes have also lower values of $p(G,k)$ than the average of randomly generated (regular) graphs. However, it is a little surprising that the Harary graphs have not  better performance even than the average among all the randomly generated regular graphs.
\begin{figure}[htbp]
\begin{center}
\includegraphics[height=8cm]{3}\\
{\bf Figure 3.}\,\ The values of $p(G,k)$ for $Q_6,\overline{R}_6,\overline{R}R_6$, and $H(4^3,3)$.
\end{center}
\end{figure}

\section{Conclusion}
Graphs with logarithmical average degree (LAD) form an important graph class, which exhibits many interesting properties in both theoretical and practical sides. We considered two probability models $p(G,k)$ and $t(G,p)$ for this graph class and chose  Harary graphs and Hypercubes as two instances, in an attempt to reveal some further properties of the class. We deduced an explicit formula of $p(G,k)$ for  Harary graphs and proved that the phase transition on $p(G,k)$ occurs in the class of Harary graphs with  logarithmical degree, which would be the first known graph class that has such phenomenon. We also performed simulations on  Hypercubes and randomly generated graphs. The numerical results suggest that the regularity is apparently a key structural feature that guarantees lower $p(G,k)$. As an extensively  studied and widely applied graphs, the Hypercubes have also better behaviour on  $p(G,k)$ than the average of the regular graphs of the same scale. In contrast, though the Harary graphs has nice properties related to graph connectedness, it is a little surprising that Harary graph, in the sense of the two measurements $p(G,k)$ and $t(G,p)$,  does not appear better behaviour than even the average of the regular graphs.

In fact, for $p(G,k)$ and some special $k$, we can observe that $p(H(4^w,w),k)>p(Q_n,k)$. For a $d$-regular graph $G$ of order $N$, it is clear that $p(G,N-3)={N\choose 3}-N{d\choose 2}+2c_3(G)$, where $c_3(G)$ is the number of cycles of length 3 in $G$. Therefore, $p(H(4^w,w),4^w-3)>p(Q_{2w},4^w-3)$ as $c_3(H(4^w,w))>c_3(Q_{2w})=0$. This observation can be also extended in terms of graph girth (the length of a shortest cycle). Even so, comparing the value $p(G,k)$ of two graphs for general $k$ is very difficult in general.

Finally, motivated by  the phase transition in Harary graphs, we propose the following two problems on the existence of the phase transition of $p(G,k)$ for Hypercubes and regular graphs:
\begin{problem}Is there a real number $c$ with $0<c<1$ such that $p(Q_n,k)=0$ when $k<c2^n$ and $p(Q_n,k)=1$ when $k>c2^n$, as $n\rightarrow\infty$?
\end{problem}
\begin{problem}Let ${\cal G}_n$ be the class of all the $n$-regular connected graphs with $2^n$ vertices. Is there a real number $c$ with $0<c<1$ such that, for $G_n\in{\cal G}_n$,  $p(G_n,k)=0$ when $k<c2^n$ and $p(G_n,k)=1$ when $k>c2^n$, almost certainly as $n\rightarrow\infty$?
\end{problem}
We remark that, by our theoretical results for Harary graphs, and simulation results for Hypercubes and randomly generated regular graphs,  if the $c$ in the two problems above exists, then it should satisfy $c>0.25$.

\section*{Declaration of competing interest}
 The authors declare that they have no known competing financial interests or personal relationships that could have appeared to influence the work reported in this paper.
 \section*{Acknowledgements}
This work was supported by the National Natural Science Foundation of China [Grant number, 12361070].

\section{Appendix}
\setlength{\tabcolsep}{1.8pt}
\footnotesize
\begin{longtable}{l *{12}{c}}
\caption{The numerical values of $\lambda_k(G)/4000$}
\label{tab:lambda}
\endfirsthead
\endhead
\endfoot
\endlastfoot

& \multicolumn{12}{c}{$k=1\sim12$} \\
\toprule
\diagbox{$G$}{$k$} & 1 & 2 & 3 & 4 & 5 & 6 & 7 & 8 & 9 & 10 & 11 & 12 \\
\midrule
$Q_6$ & 0 & 0 & 0 & 0 & 0 & 0 & 0 & 0 & 0 & 0 & 0 & 0.0005 \\
$\overline{R}_6$ & 0.0156 & 0.0238 & 0.0363 & 0.0498 & 0.0655 & 0.0805 & 0.0983 & 0.1163 & 0.1350 & 0.1545 & 0.1753 & 0.1965 \\
$\overline{R}R_6$ & 0 & 0.0010 & 0.0015 & 0.0025 & 0.0043 & 0.0060 & 0.0085 & 0.0118 & 0.0160 & 0.0215 & 0.0283 & 0.0368 \\
$H(64,3)$ & 0 & 0 & 0 & 0 & 0 & 0 & 0.0003 & 0.0005 & 0.0015 & 0.0045 & 0.0105 & 0.0163 \\

\bottomrule[0.6pt]
\\[2pt]

& \multicolumn{12}{c}{$k=13\sim24$} \\
\toprule
\diagbox{$G$}{$k$} & 13 & 14 & 15 & 16 & 17 & 18 & 19 & 20 & 21 & 22 & 23 & 24 \\
\midrule
$Q_6$ & 0.0010 & 0.0015 & 0.0035 & 0.0063 & 0.0088 & 0.0115 & 0.0153 & 0.0235 & 0.0303 & 0.0445 & 0.0605 & 0.0770 \\
$\overline{R}_6$ & 0.2200 & 0.2443 & 0.2688 & 0.2950 & 0.3233 & 0.3523 & 0.3803 & 0.4100 & 0.4405 & 0.4710 & 0.5025 & 0.5335 \\
$\overline{R}R_6$ & 0.0470 & 0.0603 & 0.0763 & 0.0948 & 0.1183 & 0.1433 & 0.1748 & 0.2105 & 0.2505 & 0.2965 & 0.3450 & 0.3995 \\
$H(64,3)$ & 0.0303 & 0.0503 & 0.0745 & 0.1155 & 0.1623 & 0.2178 & 0.2870 & 0.3695 & 0.4443 & 0.5363 & 0.6160 & 0.6985 \\
\bottomrule[0.6pt]
\\[2pt]

& \multicolumn{12}{c}{$k=25\sim36$} \\
\toprule
\diagbox{$G$}{$k$} & 25 & 26 & 27 & 28 & 29 & 30 & 31 & 32 & 33 & 34 & 35 & 36 \\
\midrule
$Q_6$ & 0.0975 & 0.1205 & 0.1503 & 0.1823 & 0.2208 & 0.2620 & 0.3100 & 0.3643 & 0.4188 & 0.4790 & 0.5360 & 0.5953 \\
$\overline{R}_6$ & 0.5640 & 0.5953 & 0.6258 & 0.6568 & 0.6863 & 0.7158 & 0.7435 & 0.7708 & 0.7958 & 0.8198 & 0.8423 & 0.8628 \\
$\overline{R}R_6$ & 0.4568 & 0.5158 & 0.5745 & 0.6350 & 0.6940 & 0.7475 & 0.7980 & 0.8420 & 0.8805 & 0.9133 & 0.9383 & 0.9573 \\
$H(64,3)$ & 0.7728 & 0.8298 & 0.8845 & 0.9235 & 0.9550 & 0.9713 & 0.9835 & 0.9915 & 0.9960 & 0.9978 & 0.9985 & 0.9995 \\
\bottomrule[0.6pt]
\\[2pt]

& \multicolumn{12}{c}{$k=37\sim48$} \\
\toprule
\diagbox{$G$}{$k$} & 37 & 38 & 39 & 40 & 41 & 42 & 43 & 44 & 45 & 46 & 47 & 48 \\
\midrule
$Q_6$ & 0.6608 & 0.7170 & 0.7728 & 0.8230 & 0.8658 & 0.9078 & 0.9293 & 0.9555 & 0.9688 & 0.9788 & 0.9845 & 0.9888 \\
$\overline{R}_6$ & 0.8810 & 0.8985 & 0.9140 & 0.9280 & 0.9400 & 0.9503 & 0.9595 & 0.9668 & 0.9733 & 0.9785 & 0.9825 & 0.9860 \\
$\overline{R}R_6$ & 0.9998 & 0.9998 & 1 & 1 & 0.9998 & 0.9998 & 0.9998 & 0.9998 & 0.9995 & 0.9990 & 0.9978 & 0.9943 \\
$H(64,3)$ & 1 & 1 & 1 & 1 & 1 & 1 & 1 & 1 & 1 & 1 & 1 & 1 \\
\bottomrule[0.6pt]
\\[2pt]

& \multicolumn{12}{c}{$k=49\sim60$} \\
\toprule
\diagbox{$G$}{$k$} & 49 & 50 & 51 & 52 & 53 & 54 & 55 & 56 & 57 & 58 & 59 & 60 \\
\midrule
$Q_6$ & 0.9945 & 0.9958 & 0.9970 & 0.9975 & 0.9983 & 0.9990 & 0.9993 & 0.9990 & 0.9990 & 0.9980 & 0.9963 & 0.9925 \\
$\overline{R}_6$ & 0.9885 & 0.9910 & 0.9925 & 0.9935 & 0.9948 & 0.9953 & 0.9960 & 0.9963 & 0.9958 & 0.9948 & 0.9933 & 0.9883 \\
$\overline{R}R_6$ & 0.9998 & 0.9998 & 1 & 1 & 0.9998 & 0.9998 & 0.9998 & 0.9998 & 0.9995 & 0.9990 & 0.9978 & 0.9943 \\
$H(64,3)$ & 1 & 1 & 1 & 1 & 1 & 1 & 1 & 1 & 1 & 0.9998 & 0.9993 & 0.9975 \\
\bottomrule[0.6pt]
\\[2pt]

& \multicolumn{12}{c}{$k=61\sim63$} \\
\multicolumn{13}{c}{
\begin{tabular*}{\linewidth}{@{\extracolsep{\fill}} l *{4}{c} @{}}
\toprule
\diagbox{$G$}{$k$} & 61 & 62 & 63 \\
\midrule
$Q_6$ & 0.9775 & 0.9048 & 1 \\
$\overline{R}_6$ & 0.9725 & 0.9048 & 1 \\
$\overline{R}R_6$ & 0.9810 & 0.9048 & 1\\
$H(64,3)$ & 0.9848 & 0.9048 & 1 \\
\bottomrule[0.6pt]
\end{tabular*}
}
\\[2pt]
\end{longtable}

\end{document}